\documentclass[10pt]{amsart}
\usepackage{amsfonts,latexsym, color}
\usepackage{amsmath}
\usepackage{amssymb}

 \font\tenmsb=msbm10 at 12pt \font\sevenmsb=msbm7 at 8pt \font\fivemsb=msbm5 at
6pt
\newfam\msbfam
\textfont\msbfam=\tenmsb \scriptfont\msbfam=\sevenmsb \scriptscriptfont\msbfam=\fivemsb

\def\R{{\mathbb R}}

\def\B{{\mathbb B}}

\def\L{{\mathcal L}}

\begin{document}
\newcommand{\reset}{\setcounter{equation}{0}}

\newcommand{\beq}{\begin{equation}}
\newcommand{\noi}{\noindent}
\newcommand{\eeq}{\end{equation}}
\newcommand{\dis}{\displaystyle}
\newcommand{\mint}{-\!\!\!\!\!\!\int}

\def \theequation{\arabic{section}.\arabic{equation}}

\newtheorem{thm}{Theorem}[section]
\newtheorem{lem}[thm]{Lemma}
\newtheorem{cor}[thm]{Corollary}
\newtheorem{prop}[thm]{Proposition}
\theoremstyle{definition}
\newtheorem{defn}[thm]{Definition}
\newtheorem{rem}[thm]{Remark}

\def \bx{\hspace{2.5mm}\rule{2.5mm}{2.5mm}} \def \vs{\vspace*{0.2cm}} \def
\hs{\hspace*{0.6cm}}
\def \ds{\displaystyle}
\def \p{\partial}
\def \O{\Omega}
\def \H{{\mathbb H}}
\def \b{\beta}
\def \m{\mu}
\def \T{{\mathbb T}}
\def \ou{{\overline u}}
\def \ov{{\overline v}}
\def \D{\Delta}
\def \M{{\mathcal M}}
\def \l{\lambda}
\def \s{\sigma}
\def \e{\varepsilon}
\def \a{\alpha}
\def \o{\omega}
\def \b{\beta}
\def \wv{{\widetilde v}}
\def \hw{{\widehat w}}
\def \HH{{\mathcal H}}
\def \E{{\mathbb E}}
\def\cqfd{%
\mbox{ }%
\nolinebreak%
\hfill%
\rule{2mm} {2mm}%
\medbreak%
\par%
}
\def \pr {\noindent {\bf Proof.} }
\def \rmk {\noindent {\it Remark} }
\def \esp {\hspace{4mm}}
\def \dsp {\hspace{2mm}}
\def \ssp {\hspace{1mm}}

\title{First order Hardy inequalities revisited}
\author{Xia Huang}
\address{School of Mathematical Sciences, East China Normal
University and Shanghai Key Laboratory of PMMP, Shanghai, 200241, China}
\email{xhuang@cpde.ecnu.edu.cn}
\author{Dong Ye}
\address{School of Mathematical Sciences, East China Normal
University and Shanghai Key Laboratory of PMMP, Shanghai, 200241, China}
\address{IECL, UMR 7502, University of Lorraine, 57050 Metz, France}
\email{dye@math.ecnu.edu.cn}

\date{}

\begin{abstract}
In this paper, we consider the first order Hardy inequalities using simple equalities. This basic setting not only permits to derive quickly many well-known Hardy inequalities with optimal constants, but also supplies improved or new estimates in miscellaneous situations, such as multipolar potential, the exponential weight, hyperbolic space, Heisenberg group, the edge Laplacian, or the Grushin type operator.
\end{abstract}


\maketitle
\begin{center}
{\sl Dedicated to Professor Chaojiang Xu's 65th birthday}
\end{center}
\section{Introduction}
The Hardy inequalities go back to G.H.~Hardy, who showed in \cite{ha1} a very famous estimate: Let $p > 1$, then
\begin{align}
\label{Hardy0}
\int_0^\infty |u'(x)|^p dx \geq \left(\frac{p-1}{p}\right)^p \int_0^\infty \frac{|u(x)|^p}{x^p} dx, \quad \forall \; u \in  C^1(\R_+), \; u(0) = 0.
\end{align}
Since one century, the Hardy type inequalities have been enriched extensively and broadly, they play important roles in many branches of analysis and geometry. More generally, we call first order Hardy inequalities, the estimates like
\begin{align*}
\int_{\M} V|\nabla u|^p d\mu \geq \int_{\M} W|u|^p d\mu,
\end{align*}
with positive weights $V$, $W$ and $u$ in suitable function spaces. They are also called weighted Poincar\'e or Hardy-Poincar\'e inequalities.

\medskip
A huge literature exists on the Hardy inequalities, it's just impossible nor our intension to mention all the progress even for the first order case, we refer to the classical and recent books \cite{ma1, ok, bel, gm2, rs} for interested readers. The modest objective here is to show that many first order Hardy inequalities can be derived naturally and quickly from simple equalities, which on one hand yield classical Hardy inequalities with optimal constants; and on the other hand provide new Hardy inequalities or improve some well-known results.

\medskip
The most well-known first order Hardy inequality is the following: Let $\alpha \in \R$,
\begin{align}
\label{cHR1a}
\int_{\R^n} \frac{|\nabla u|^2}{|x|^\alpha} dx \geq \frac{(n-2 - \alpha)^2}{4}\int_{\R^n}\frac{u^2}{|x|^{\alpha +2}} dx, \quad \forall\; u \in C_c^1(\R^n\backslash\{0\}).
\end{align}
Here and after, $|\cdot|$ denotes the Euclidean norm. The optimal constant could be firstly shown in \cite[page 259]{hlp} with $n=1$. For simplicity, we consider only real valued functions; and without special remark, the functions $u$ are $C^1$, compactly supported away from the singularities of involved weights. In general, applying density argument, many estimates hold true in larger functional spaces.

\medskip
The inequality \eqref{cHR1a} can be seen as a direct consequence of the following equality (see for instance \cite{dv} with $\a = 0$ and \cite[Lemma 2.3(i)] {tz} with $\alpha = 2$). For any $u \in C_c^1(\R^n\backslash\{0\})$ and $\alpha \in \R$, if $v = |x|^\frac{n-2-\alpha}{2}u$, there holds
\begin{align}
\label{HR1a}
\int_{\R^n}\frac{|\nabla u|^2}{|x|^\alpha} dx = \frac{(n-2-\alpha)^2}{4} \int_{\R^n}\frac{u^2}{|x|^{\alpha+2}} dx + \int_{\R^n}\frac{|\nabla v|^2}{|x|^{n-2}} dx.
\end{align}
There is also a radial derivative version of the above equality: $\forall\; u \in C_c^1(\R^n\backslash\{0\})$,
\begin{align}
\label{rHR1a}
\int_{\R^n}\frac{|\p_r u|^2}{|x|^\alpha} dx = \frac{(n-2 -\alpha)^2}{4} \int_{\R^n}\frac{u^2}{|x|^{\alpha+2}} dx + \int_{\R^n}\frac{|T_\a(u)|^2}{|x|^{\a}} dx,
\end{align}
where $\p_r$ is the radial derivative $\p_r = \frac{x\cdot\nabla}{r}$, and $T_\a = \p_r + \frac{n-2-\alpha}{2r}$.

\medskip
Indeed, the equalities \eqref{HR1a} and \eqref{rHR1a} are equivalent, since
$$|\nabla u|^2 = |\p_r u|^2 + \sum_{j = 1}^n |L_j u|^2, \quad \mbox{where }\; L_j u = \p_ju - \frac{x_j}{r}\p_r u, \;\; \forall\; 1\leq j \leq n.$$
Therefore
$$\left|\nabla\left(|x|^\frac{n-2-\alpha}{2}u\right)\right|^2 = |x|^{n-2 - \a}\Big(|T_\a(u)|^2 + \sum_{j = 1}^n |L_j u|^2\Big).$$
From \eqref{HR1a} and \eqref{rHR1a}, as the remainder terms are nonnegative, we deduce easily the optimal Hardy inequalities \eqref{cHR1a}, and the non-attainability of equality in \eqref{cHR1a} over the corresponding Banach space, for example,  $D^{1,2}(\R^n)$ if $\a = 0$ and $n \geq 3$.

\medskip
The formula \eqref{rHR1a} was found firstly by Brezis-V\'azquez in 1997 for $\a = 0$. They called it {\it the ``magical" computation} in \cite[page 454]{bv}, and applied it to get a famous improved Hardy inequality: For any bounded domain $\O \subset \R^n$ with $n \geq 2$,
\begin{align}
\label{BV}
\int_\O |\nabla u|^2 dx - \frac{(n-2)^2}{4}\int_\O \frac{|u|^2}{|x|^2} dx \geq H_2\left(\frac{\omega_n}{|\O|}\right)^\frac{2}{n}\int_\O |u|^2 dx, \quad \forall\; u \in H_0^1(\O).
\end{align}
Here $H_2$ stands for the first eigenvalue of $-\D$ in $H_0^1(\B^2)$. In this paper, $\B^n$ denotes always the Euclidean unit ball in $\R^n$. The above inequality is the departure point for a rich literature on the improvement of classical Hardy inequalities, over bounded or not domains in $\R^n$.

\medskip
More recently, in the important works of Ghoussoub-Moradifam \cite{gm1, gm2}, they introduced a notion of Bessel pair to study Hardy type inequalities. Let $V$, $W$ be positive $C^1$ radial functions on $\B^n\backslash \{0\}$ such that
\begin{align}
\label{B1}
\int_0^1 \frac{1}{r^{n-1}V(r)}dr + \int_0^1 r^{n-1}V(r) dr < \infty,
\end{align}
they proved that the following assertions are equivalent:
\begin{itemize}
\item[(i)] $(V, W)$ is a $n$-dimensional Bessel pair on $(0, 1)$ with constant $c> 0$, that is, the ODE
\begin{align}
\label{B2}
y''(r) + \left(\frac{n-1}{r}+ \frac{V'}{V}\right)y'(r) + \frac{cW}{V}y(r) = 0
\end{align}
has a positive solution on the interval $(0, 1)$.
\item[(ii)] For all $u \in C_c^1(\B^n)$, there holds
\begin{align}
\label{B3}
\int_{\B^n} V|\nabla u|^2 dx \geq \b_{V,W}\int_{\B^n} Wu^2 dx.
\end{align}
Here $\b_{V,W}$ is the best constant, given as the supermum of $c$ satisfying (i).
\end{itemize}

Ghoussoub-Moradifam used the above Bessel pair notion to improve, extend and unify many results on weighted Hardy inequalities, and studied the corresponding best constants. Their ideas have inspired many researches in the last decade and turn out to be very successful.

\medskip
Here we attempt a more general approach using another point of view, where the departure point is an easy equality as follows. Let $(\mathcal M, g)$ be a Riemannian manifold and $\Omega\subset \mathcal M$ be open. Consider $V\in C^1(\Omega)$ and $\vec{F}\in C^1(\Omega, T_g{\mathcal M})$, then for any $u\in C_c^1(\Omega)$, any family of inner product $\langle \cdot, \cdot\rangle \in C^1(\O, \Lambda^2 T_g{\mathcal M})$ (not necessarily that given by $g$), there holds
\begin{align}
\label{E1}
\int_{\Omega} V \|\nabla_g u\|^2 dg = \int_{\Omega} \left[{\rm div}_g(V \vec{F})-V\|\vec {F}\|^2\right]u^2 dg + \int_{\Omega} V\left\|\nabla_g u + u\vec {F}\right\|^2 dg.
\end{align}
The formula \eqref{E1} is proved just by developing $\big\|\nabla_g u + u\vec {F}\big\|^2$ and integration by parts. Here $\nabla_g$ and $dg$ are respectively the gradient and volume form with respect to the metric $g$; $\|\cdot\|$ denotes the norm associated to $\langle \cdot, \cdot\rangle$; and ${\rm div}_g$ is the adjoint operator of $\nabla_g$ with respect to $\langle \cdot, \cdot\rangle$ and $g$.

\medskip
In particular, let $\vec F=-\frac{\nabla_g f}{f}$ with $f$ positive in $C^2(\Omega)$, we see that
\begin{align}
\label{E2}
\int_{\Omega} V \|\nabla_g u\|^2 dg  = - \int_{\Omega} \frac{{\rm div}_g(V \nabla_g f)}{f}u^2 dg + \int_{\Omega} V\left\|\nabla_g u - \frac{u}{f}\nabla_g f \right\|^2 dg.
\end{align}
If moreover $V = 1$, and the inner product $\langle \cdot, \cdot\rangle$ is that given by $g$, we arrive at
\begin{align}
\label{E3}
\int_{\Omega} \|\nabla_g u\|^2 dg  = - \int_{\Omega} \frac{\D_g f}{f}u^2 dg + \int_{\Omega} f^2\left\|\nabla_g \left( \frac{u}{f}\right) \right\|^2 dg.
\end{align}

\medskip
The idea to use \eqref{E1}, \eqref{E2} or \eqref{E3} to get Hardy inequalities has been used in various situations, for example Lemma 2.1 in \cite{ftt} presents \eqref{E1} in Euclidean spaces. The following are some classical or more recent Hardy inequalities which can be derived by the equalities \eqref{E1}--\eqref{E3}.

\begin{itemize}
\item Taking $f = |x|^\frac{2-n+\alpha}{2}$ and $V = |x|^{-\a}$ in \eqref{E2}, we get $W = -\frac{{\rm div}(V \nabla f)}{f} = \frac{(n-2-\alpha)^2}{4|x|^{\a+2}}$ in $\R^n\backslash\{0\}$ which gives \eqref{HR1a}. In other words, equality \eqref{HR1a} is a special example of \eqref{E2}.
\item In 1934, J. Leray \cite{le} proved a famous inequality in dimension two:
\begin{align}
\label{Leray}
\int_{\B^2} |\nabla u|^2 dx \geq \frac{1}{4} \int_{\B^2}\frac{u^2}{|x|^2(\ln |x|)^2} dx, \quad \forall\; u \in H^1_0(\B^2).
\end{align}
Wang-Willem showed in \cite{ww} a generalization of Leray's inequality: Let $n \geq 1$, $\alpha \leq n - 2$, $\Omega=\R^n\backslash \overline{\B^n}$ or $\B^n$. Then for any $u \in C^1_c(\Omega)$,
\begin{align}\label{HWW}
\int_\Omega \frac{|\nabla u|^2}{|x|^\alpha} dx \geq \frac{(n-2 - \alpha)^2}{4} \int_\Omega\frac{u^2}{|x|^{\alpha+2}} dx + \frac{1}{4}\int_\Omega\frac{u^2}{|x|^{\alpha+2}(\ln|x|)^2} dx.
\end{align}
The estimate \eqref{HWW} can be proved using $f = |x|^\frac{2-n+\a}{2}\left|\ln|x|\right|^\frac{1}{2}$, $V = |x|^{-\a}$ in \eqref{E2}.
\item Carron \cite {c} proved the following well-known Hardy inequality for non compact manifold: Let $(\M, g)$ be a non parabolic complete Riemannian manifold such that there is a nonnegative function $\rho$ satisfying $\|\nabla \rho\| = 1$; $\Delta_g \rho \geq -C/\rho$ in ${\mathcal D}'(\M)$ with $C > 1$;  and the zero set of $\rho$ is compact with null capacity. Then we have
$$\int_{\M} \|\nabla u\|^2 dx \geq \frac{(C-1)^2}{4} \int_{\M}\frac{u^2}{\rho^2} dx, \quad \forall\; u \in C_c^1(\M).$$
His proof is equivalent to use $f(\rho) = \rho^{\frac{1-C}{2}}$ in \eqref{E3}, hence $W = -\frac{\Delta_g f}{f} \geq \frac{(C-1)^2}{4\rho^2}$.
\item Brock-Chiacchio-Mercaldo \cite{bcm} used the following inequality to get some weighted isoperimetric inequalities in cones: Let $n \geq 3$, $\O = \{x_1 > 0, x_n > 0\}$, $d\mu = x_1^k|x|^mdx$, $k \geq 1$ and $m\in \R$. Denote $E$ the closure of $C^1(\overline\O)$ functions vanishing on $\{x_1 = 0\}$ in $H^1(\O, d\mu)$, then
$$\int_\O |\nabla u|^2 d\mu \geq \left[\frac{(n+m+k)^2}{4} - m\right]\int_\O \frac{u^2}{|x|^2} d\mu, \quad \forall\; u \in E.$$
This can be derived from \eqref{E2} by taking $V = x_1^k|x|^m$ and $f(x) = x_1|x|^{-\frac{n+m+k}{2}}$.
\item Let $f(x) = y(r)$, we can check that the equation \eqref{B2} is equivalent to
$$-{\rm div}(V\nabla f) = cWf \quad \mbox{in }\; \B^n\backslash\{0\}.$$
By \eqref{E2}, for any $u \in C_c^1(\B^n\backslash\{0\})$, the Hardy inequality \eqref{B3} holds true with $c$ tending to $\b_{V,W}$. The technical assumption \eqref{B1} permits to extend \eqref{B3} for $u \in H^1_0(\B^n)$. Conversely, as indicated by \cite[Lemma 2.4] {gm1}, the validity of \eqref{B3} yields the existence of a positive function $f$ such that $-{\rm div}(V\nabla f)\geq cWf$ in $\B^n\backslash\{0\}$.
\end{itemize}

As illustrated by above examples, the equalities \eqref{E1}--\eqref{E3} permit to find very quickly many first order Hardy inequalities. More examples will be displayed later. The best constants can be obtained with suitable choice of $\vec{F}$, $f$, $V$ and optimization on the sequel weight $W$. For example, to obtain \eqref{HWW}, we test $f(x) = |x|^\beta\left|\ln|x|\right|^\gamma$ with $\beta, \gamma \in \R$ and $V = |x|^{-\a}$. There holds, for any $|x|\ne 1$,
$$
-\frac{{\rm div}(V\nabla f)}{f}= -\frac{\beta (n -2 - \a +\beta)}{|x|^{2+\a}} - \frac{\gamma (n -2 -\alpha + 2\beta)}{|x|^{2+\a}\left|\ln|x|\right|} + \frac{\gamma(1-\gamma)}{|x|^{2+\a}(\ln|x|)^2}.
$$
To maximize the coefficient for the first term on r.h.s., we choose $\beta = \frac{2 - n +\a}{2}$, which cancels the second term, and we optimize the last one with $\gamma = \frac{1}{2}$. Hence \eqref{HWW} is valid for $u \in C_c^1(\O\backslash\{0\})$, we conclude eventually with the density argument if $\O = \B^n$.

\medskip
Furthermore, in \eqref{E2} or \eqref{E3}, we remark that the equality holds true without the last term, if and only if $u/f$ is locally constant, this resolves quickly the attainability issue, and suggests us the idea to construct appropriate approximation sequence showing the eventual optimality of involved constants.

\medskip
The use of vector fields $\vec{F}$ in \eqref{E1} can be helpful to handle anisotropic situation (even the final choice could be a gradient field). For example, Tidblom \cite{ti} and Filippas-Tertikas-Tidblom \cite{ftt} used \eqref{E1} to get various Hardy-inequalities on orthogonal cones in $\R^n$. See also the consideration for multipolar Hardy inequalities in \cite{bde, cz}.

\medskip
In the sequel, we will handle many other situations just using the identity \eqref{E1} (hence \eqref{E2} or \eqref{E3}), and we show miscellaneous improved or new Hardy inequalities.
\begin{itemize}
\item For multipolar potential case, we improve and generalize several estimates given by Felli-Marchini-Terracini [17], Bosi-Dolbeault-Esteban \cite{bde}, Cazacu-Zuazua \cite{cz};
\item We answer an open problem of Blanchet-Bonforte-Dolbeaul-Grillo-V\'azquez \cite{bbdgv} (see also \cite{gm1}), by showing best constants for Hardy inequality with weight $(1+|x|^2)^\beta$;
\item For the exponential weight in $\R^n$, we improve an inequality of Escobar-Kavian \cite{ek};
\item We show an improved Hardy inequality on the unit disc of hyperbolic space;
\item We consider also Heisenberg group, the Grushin type operators or the edge Laplacian, and supply various examples of new Hardy inequalities.
    \end{itemize}
At last, we give some discussion for more general $L^p$ setting and Bessel pair.

\section{Multipolar Hardy inequalities}
\reset
Consider $n \geq 3$ and $m$ distinct points $(a_i)_{1\leq i \leq m}$ in $\R^n$, if we denote $r_i = |x - a_i|$, Felli-Marchini-Terracini studied the multipolar Hardy inequality with $\l_i > 0$ as follows.
\begin{align}
\label{HFMT}
\int_{\R^n}|\nabla u|^2 dx \geq \sum_{i = 1}^m \l_i\int_{\R^n} \frac{u^2}{r_i^2}dx, \quad \forall\; u \in H^1(\R^n).
\end{align}
They proved that \eqref{HFMT} holds true if and only if $\l_1 +\ldots + \l_m \leq \frac{(n-2)^2}{4},$ see \cite[Theorem 1.1]{fmt}.

\medskip
In the special borderline case when all $\l_i$ are equal to $\frac{(n-2)^2}{4m}$, Bosi-Dolbeault-Esteban improves \eqref{HFMT} in \cite[Lemma 8]{bde}: For any $u \in H^1(\R^n)$,
\begin{align}
\label{mHR2}
\int_{\R^n}|\nabla u|^2 dx \geq \frac{(n-2)^2}{4m} \left[\sum_{i = 1}^m\int_{\R^n} \frac{u^2}{r_i^2}dx + \frac{1}{m}\sum_{1\leq i < j \leq m} \int_{\R^n}\frac{|a_i-a_j|^2}{r_i^2r_j^2}u^2 dx\right].
\end{align}

Here we present a multipolar Hardy inequality which improves and generalizes completely \eqref{HFMT} and \eqref{mHR2}.
\begin{thm}
Let $n \geq 3$, and $(\l_i) \in (0, \infty)^m$ satisfy
\begin{align}
\label{mHR0}
\sum_{i = 1}^m \l_i = \frac{(n-2)^2}{4}.
\end{align}
Then for any $u \in H^1(\R^n)$, and any family of $m$ distinct points $(a_i)_{1\leq i \leq m}$ in $\R^n$, there holds
\begin{align}
\label{mHR1}
\int_{\R^n}|\nabla u|^2 dx \geq \sum_{i = 1}^m \l_i\int_{\R^n} \frac{u^2}{r_i^2}dx + \frac{4}{(n-2)^2}\sum_{1\leq i < j \leq m}\l_i\l_j\int_{\R^n}\frac{|a_i-a_j|^2}{r_i^2r_j^2}u^2 dx,
\end{align}
where $r_i = |x-a_i|$.
\end{thm}

\medskip
\pr Let $\alpha = (\a_i) \in \R^m$ and
$$\vec F_\alpha(x) = \sum_{1\leq i \leq m}\alpha_i\frac{x-a_i}{r_i^2}.$$
Then direct calculation gives
\begin{align*}
{\rm div}(\vec F_\a) - |\vec F_\a|^2 & = \sum_{1\leq i \leq m}\frac{\alpha_i(n-2)}{r_i^2} - \sum_{1\leq i \leq m}\frac{\alpha_i^2}{r_i^2} - \sum_{1\leq i \ne j \leq m}\alpha_i\alpha_j \frac{\langle x-a_i, x-a_j\rangle}{r_i^2r_j^2}\\
& = \sum_{1\leq i \leq m}\frac{\alpha_i(n-2) -\a_i^2}{r_i^2} + \frac{1}{2}\sum_{1\leq i \ne j \leq m}\left(\frac{|a_i-a_j|^2}{r_i^2r_j^2} - \frac{1}{r_i^2} - \frac{1}{r_j^2}\right)\a_i\a_j\\
& = \sum_{1\leq i \leq m}\frac{\alpha_i(n-2 -S_\a)}{r_i^2} + \frac{1}{2}\sum_{1\leq i \ne j \leq m}\frac{|a_i-a_j|^2}{r_i^2r_j^2} \a_i\a_j.
\end{align*}
Here $S_\a = \sum_{1\leq i \leq m} \a_i$, and we used the equality
$$\frac{|a_i-a_j|^2}{r_i^2r_j^2} = \frac{|(x-a_i)-(x-a_j)|^2}{r_i^2r_j^2} = \frac{1}{r_i^2} + \frac{1}{r_j^2}  -2\frac{\langle x-a_i, x-a_j\rangle}{r_i^2r_j^2}.$$

Now we want to choose $\l_i = (n-2 -S_\a)\a_i$ with $\alpha_i > 0$ and $S_\a < n-2$. By \eqref{mHR0},
$$\frac{(n-2)^2}{4} = (n-2 -S_\a)S_\a, \quad \mbox{hence }\; S_\alpha = \frac{n-2}{2}, \;\; \alpha_i = \frac{2\l_i}{n-2}.$$
Finally, with the above values of $(\a_i)$, we obtain
\begin{align*}
{\rm div}(\vec F_\a) - |\vec F_\a|^2 = \sum_{1\leq i \leq m}\frac{\l_i}{r_i^2} + \frac{4}{(n-2)^2}\sum_{1\leq i < j \leq m}\l_i\l_j\frac{|a_i-a_j|^2}{r_i^2r_j^2} .
\end{align*}
By \eqref{E1}, the estimate \eqref{mHR1} holds true for $u \in C_c^1(\R^n\backslash\{a_i, 1\leq i \leq m\})$, and remains valid in $H^1(\R^n)$ by capacity or approximation argument as $n \geq 3$.\qed

\medskip
We should mention that the method in \cite{bde} was somehow the vector field approach as \eqref{E1}, they used the special $\vec F_\a$ where all $\a_i = \frac{n-2}{2m}$. Recently, using $\vec F_\a$ with all $\a_i = \frac{n-2}{m}$, that is
$$\vec F(x) = \frac{n-2}{m}\sum_{1\leq i \leq m}\frac{x-a_i}{r_i^2},$$
Cazacu-Zuazua obtained in \cite{cz} another optimal multipolar Hardy inequality:
\begin{align}
\label{mHR3}
\int_{\R^n}|\nabla u|^2 dx \geq \frac{(n-2)^2}{m^2} \sum_{1\leq i < j \leq m} \int_{\R^n}\frac{|a_i-a_j|^2}{r_i^2r_j^2}u^2 dx, \quad \forall\; u \in H^1(\R^n).
\end{align}
The following is a family of Hardy inequalities which generalizes \eqref{mHR3}.
\begin{thm}
Let $n \geq 3$, and $ (\mu_i) \in (0, \infty)^m$ satisfy $\sum_{1\leq i \leq m} \mu_i = n-2$. Then for any $u \in H^1(\R^n)$ and any family of $m$ distinct points $(a_i)_{1\leq i \leq m}$ in $\R^n$, there holds
\begin{align}
\label{mHR4}
\int_{\R^n}|\nabla u|^2 dx \geq \sum_{1\leq i < j \leq m}\mu_i\mu_j\int_{\R^n}\frac{|a_i-a_j|^2}{r_i^2r_j^2}u^2 dx,
\end{align}
where $r_i = |x-a_i|$.
\end{thm}

The proof uses exactly the same computation as for \eqref{mHR1}, with now $\alpha = (\mu_i)$ and $S_\a = n-2$, so we omit it.

\section{Best constant with weight $(1 + |x|^2)^\beta$}
\reset
To study the asymptotic behavior of solutions to fast diffusion equation via entropy estimate \cite{bbdgv, bbdgv2} Blanchet-Bonforte-Dolbeaul-Grillo-V\'azquez showed some spectral gap estimates, using Hardy inequalities with weights $(1 + |x|^2)^\beta$. More precisely, they studied the following Hardy type estimate: Let $n \geq 3$, $\a \in \R$,
\begin{align}
\label{HBBDGV}
\int_{\mathbb{R}^n} (1+|x|^2)^\a |\nabla u|^2 dx \geq C\int_{\mathbb{R}^n}(1+|x|^2)^{\a-1}u^2 dx, \quad \forall\; u\in C_c^1(\mathbb{R}^n).
\end{align}
As indicated in \cite{bbdgv2}, under suitable initial conditions for the fast diffusion equation $\p_t v = \D(v^{1 + \frac{1}{\a}})$ in $\R^n$, the best constant in \eqref{HBBDGV} provides the sharp decay rate of the entropy. They supplied in \cite[Theorem 2]{bdgv} the best constant in \eqref{HBBDGV} for $\a < 0$, see also \cite{gm1}. The best constant with $\a > 0$ was left open, except for $\a = n$, see \cite[Table 1]{bbdgv} and \cite[Open problem (3)]{gm2}. Here we give a complete answer.

\begin{thm}
Let $n>2$, let $C_{\a, n}$ denote the best constant for \eqref{HBBDGV}. We have
\begin{align*}
C_{\a, n} = \frac{(n+2\alpha - 2)^2}{4}, \;\; \mbox{if } \; \frac{2-n}{2} \leq \a \leq \frac{n+2}{2}; \quad \mbox{and}\quad  C_{\a, n} = 2(\a - 1)n, \;\; \mbox{if } \; \a > \frac{n+2}{2}.
\end{align*}
\end{thm}

\medskip
\pr Let $V=(1+r^2)^\a$, $f=(1+r^2)^\gamma$ with $\gamma\in \mathbb{R}$ and $r =|x|$. Direct calculation yields
\begin{align*}
W := -\frac{{\rm div}(V\nabla f)}{f} & = \frac{-V'f'-V\Delta f}{f} = \left[Q(\gamma) + \frac{4\gamma(\gamma + \a -1)}{1+r^2}\right](1+r^2)^{\a - 1},
\end{align*}
where
\begin{align*}
Q(\gamma) & = -2\gamma(2\gamma + n+2\a -2)  = \frac{(n+2\a - 2)^2}{4} - \left(2\gamma + \frac{n+2\a - 2}{2}\right)^2.
\end{align*}
Hence $$\max_{\R} Q(\gamma) = \frac{(n+2\a - 2)^2}{4}, \quad \mbox{reached by } \gamma = \gamma_* = -\frac{n+2\a - 2}{4}.$$

\medskip
Fixing $\gamma = \gamma_*$, there holds
\begin{align*}
\frac{4\gamma_* (\gamma_* + \a -1)}{1+r^2} & = \frac{-(n+2\a - 2)(2\a - 2 -n)}{4(1+r^2)} = \frac{n^2 - (2\a - 2)^2}{4(1+r^2)} =: T(r).
\end{align*}
So the discussion depends on the sign of $\frac{n^2 - (2\a - 2)^2}{4}$.

\medskip
\begin{itemize}
\item If $\frac{2-n}{2} \leq \a \leq \frac{n+2}{2}$, taking $\gamma = \gamma_*$, as $\min_{\R_+}T(r) = 0$, we claim that $C_{\a, n} \geq \frac{(n+2\a - 2)^2}{4}$. Finally $C_{\a, n} = \frac{(n+2\a - 2)^2}{4}$, since it's known that $C_{\a, n} \leq \frac{(n+2\a - 2)^2}{4}$, see [20].
\item Assume now $\a > \frac{n+2}{2}$. Let $f(x) = (1+r^2)^{1 - \a}$, there holds $(1+r^2)^{1 - \a}W \equiv 2(\a - 1)n$. In the spirit of \eqref{E2}, it means that $C_{\a, n} \geq 2(\a - 1)n$. Moreover, as $(1+r^2)^{\a - 1} f^2 + (1+r^2)^{\a} |\nabla f|^2 \in L^1(\R^n)$, applying standard cut-off argument, we conclude readily that $C_{\a, n} \leq 2(\a - 1)n$, so we are done. \qed
\end{itemize}

Moreover, when $\a > \frac{n+2}{2}$, \eqref{E2} means that $\lambda(1+r^2)^{1 - \a}$ are the unique minimizers which realize the best constant in the corresponding weighted Sobolev space used by [4, 5], obtained as the completion of $C_0^\infty(\R^n)$ with respect to the norm
$$\|v\|_\a = \|v\|_{L^2(\R^n, d\mu_{\a - 1})} + \|\nabla v \|_{L^2(\R^n, d\mu_\a)}, \quad \mbox{where } \; d\mu_\beta = (1 + |x|^2)^\beta dx.$$

\begin{rem}
It is worthy to mention that when $\a < \frac{2-n}{2}$, as $(1+r^2)^{\a - 1} \in L^1(\R^n)$, the inequality \eqref{HBBDGV} is meaningful only under additional condition such as
\begin{align*}
\int_{\R^n} ud\mu_{\a-1} = \int_{\R^n} (1+r^2)^{\a - 1}u dx = 0.
\end{align*}
The best constant to the Hardy inequality \eqref{HBBDGV} under the above constrain was given for all $\a < \frac{2-n}{2}$ in \cite[Theorem 2]{bdgv}.
\end{rem}

\section{On an inequality of Maz'ya}
\reset

Maz'ya established a weighted Hardy inequality as follows, see \cite{ma2} or \cite[Theorem 1, page 214]{ma1}: For $n \geq 2$, there holds
\begin{align}
\label{HMa}
\int_{\mathbb{R}^n_+} x_n|\nabla u|^2 dx \geq \frac{1}{16}\int_{\mathbb{R}^n_+} \frac{u^2}{(x^2_{n-1}+x^2_n)^{\frac{1}{2}}}dx, \quad u \in C_c^1(\overline{\R_+^n}).
\end{align}
Consequently, consider $u(x) = |x_n|^{-\frac{1}{2}}v(x)$ with $v \in C_c^1(\R^n_+)$, we have
\begin{align}
\label{HMa2}
\int_{\R^n_+} |\nabla v|^2 dx \geq \frac{1}{4}\int_{\R^n_+} \frac{v^2}{x_n^2}dx + \frac{1}{16}\int_{\R^n_+} \frac{v^2}{x_n(x^2_{n-1}+x^2_n)^{\frac{1}{2}}}dx.
\end{align}
Later on, Tidblom raised the coefficient $\Lambda=\frac{1}{16}$ to $\frac{1}{8}$ in \eqref{HMa2}, see \cite[Corollary 3.1]{ti}. However, this does not mean that the corresponding  improvement holds for \eqref{HMa}. Because $v$ vanishes on $\p\R^n_+$ will force $u$ to be zero on $\p\R^n_+$, and the use of weight like $x_n^{-1}$ in \cite[Lemma 2.1]{ti} seems to prevent the approximation or cut-off argument with respect to $x_n$. Finally, Maz'ya-Shaposhnikova \cite{ms} showed that the optimal constant in \eqref{HMa} is equal to
\begin{align*}
\lambda := \inf_{g \in C^1([0, \pi]), \|g\|_{L^2} = 1}\int_0^\pi \left[g'(t)^2 + \frac{1}{4}g(t)^2\right]\sin t dt
\end{align*}
and $\lambda = 0.1564..$ by numerical approximation.

\medskip
Here, with \eqref{E1}, we will not reach the optimal constant, but we can improve easily \eqref{HMa} with $\frac{1}{8}$ and an explicit remainder term.
\begin{thm}
Let $n\geq 2$, then for any $u\in C_c^1(\overline{\R_+^n})$,
\begin{align}
\label{newHMa}
\begin{split}
\int_{\mathbb{R}^n_+} x_n|\nabla u|^2 dx & \geq \frac{1}{8}\int_{\mathbb{R}^n_+} \frac{u^2}{(x^2_{n-1}+x^2_n)^{\frac{1}{2}}}dx\\
& \quad + \frac{7}{32}\int_{\mathbb{R}^n_+}\left[\frac{x_n}{x^2_{n-1}+x^2_n}-\frac{x_n^2}{(x^2_{n-1}+x^2_n)^{\frac{3}{2}}}\right] u^2 dx.
\end{split}
\end{align}
\end{thm}

\medskip
\pr Consider $a, b \in \R$,
$$\vec F= \left(0,\cdot\cdot\cdot,0,\frac{b x_{n-1}}{\rho^2},\frac{a}{\rho}+\frac{ b x_n}{\rho^2}\right)$$
with $\rho=|(x_{n-1}, x_n)| > 0$ and $V=x_n$. Then for $\rho > 0$, there holds
\begin{align*}
W:={\rm div}(V \vec F)-V|\vec F|^2 & = \frac{a}{\rho} - \frac{a^2+b^2-b}{\rho^2}x_n -\frac{a(1+2b)x_n^2}{\rho^3}\\
& = \frac{a}{\rho} - \Big[a^2+b^2-b+a(1+2b)\Big]\frac{x_n}{\rho^2} + a(1+2b)\frac{x_n(\rho - x_n)}{\rho^3}.
\end{align*}
We hope that $a(1+2b)\geq 0$ and
\begin{align*}
a^2+b^2-b+a(1+2b)=\left(b+\frac{2a-1}{2}\right)^2+\frac{8a-1}{4}\leq 0.
\end{align*}
In this sense, the best choice could be $a=\frac{1}{8}$ and $b=\frac{3}{8}$, which gives
\begin{align}
\label{MaW}
W=\frac{1}{8\rho} + \frac{7}{32}\frac{x_n(\rho - x_n)}{\rho^3}.
\end{align}

Furthermore, let $u \in C_c^1(\overline{\R_+^n}\backslash\{\rho = 0\})$, even $u$ does not vanish always on $\p\R_+^n$, the equality \eqref{E1} remains true with $\Omega = \R^n_+\backslash\{\rho = 0\}$, thanks to $V = 0$ on $\p\R^n_+$. Combining with \eqref{MaW}, we see that \eqref{newHMa} is valid for any $u \in C_c^1(\overline{\R_+^n}\backslash\{\rho = 0\})$.

\medskip
Finally, let $\varphi$ be a cut-off function in $\R$, using $u_\e(x) = u(x)\big[1 - \varphi(\rho/\e)\big]$, we can prove readily that the result holds true for all $u \in C_c^1(\overline{\R_+^n})$, we omit the detail. \qed

\begin{rem}
It's interesting to notice that the vector field $\vec F$ used here is not a gradient vector field. 
\end{rem}

\section{Exponentially weighted Hardy inequalities}
\reset
The Hardy inequalities with weight $K(x) = e^\frac{|x|^2}{4}$ in $\R^n$ play important roles in the study of heat equation, see \cite{ek, vz}.
Let $V = K(x)$ and $f_\alpha(x) = |x|^\frac{2-n}{2}e^{\alpha|x|^2}$, we have
\begin{align*}
-\frac{{\rm div}(V\nabla f_\alpha)}{f_\alpha} = K(x)\left[\frac{(n-2)^2}{4|x|^2} + \frac{n-2-16\alpha}{4}  - \alpha(4\alpha+1)|x|^2\right].
\end{align*}
Hence for $n \geq 3$, $u \in H^1_K(\R^n)$ the weighted Sobolev space with respect to $K(x)dx$ as in \cite{ek, vz}, there holds
\begin{align}
\label{KHnew}
\begin{split}
\int_{\R^n} |\nabla u|^2 Kdx & \geq \frac{(n-2)^2}{4}\int_{\R^n} \frac{u^2}{|x|^2}Kdx + \frac{n-2-16\alpha}{4}\int_{\R^n} u^2Kdx\\
& \quad  -\alpha(4\alpha+1) \int_{\R^n} u^2|x|^2Kdx.
\end{split}
\end{align}

Taking $\alpha = -\frac{1}{8}$, we get
\begin{thm}
Let $n \geq 3$ and $u \in H^1_K(\R^n)$, then
\begin{align}
\label{KHnew2}
\int_{\R^n} |\nabla u|^2 Kdx \geq \frac{1}{16}\int_{\R^n} u^2|x|^2Kdx + \frac{n}{4}\int_{\R^n} u^2Kdx + \frac{(n-2)^2}{4}\int_{\R^n} \frac{u^2}{|x|^2}Kdx.
\end{align}
\end{thm}
The estimate \eqref{KHnew2} improves the consideration of Escobedo-Kavian \cite{ek} by adding the last positive term with classical potential $\frac{(n-2)^2}{4|x|^2}$.

\medskip
With $\alpha = -\frac{1}{4}$ in \eqref{KHnew}, we obtain another Hardy inequality showed by V\'azquez-Zuazua \cite[Theorem 9.1]{vz}, that is, for $u \in H^1_K(\R^n)$ with $n \geq 3$,
\begin{align}
\label{KH2}
\int_{\R^n} |\nabla u|^2 Kdx \geq \frac{(n-2)^2}{4}\int_{\R^n} \frac{u^2}{|x|^2}Kdx + \frac{n+2}{4}\int_{\R^n} u^2Kdx.
\end{align}

\section{Hyperbolic disc}
\reset
Consider the Poincar\'e hyperbolic space $\H^n$ over $\B^n$ with the metric
$$
ds^2=\frac{4dx^2}{(1- r^2)^2}, \quad r = |x|.
$$
It is well known that the volume form of $\H^n$ is $dv_\H=\frac{2^n}{(1- r^2)^n}|dx|^2$, the associated gradient operator and Laplacian are respectively
$$\nabla_\H =\frac{1-r^2}{2}\nabla, \quad \Delta_\H=\frac{1-r^2}{4}\Big[(1-r^2)\Delta+2(n-2)x\cdot\nabla \Big].$$
Recall that $\nabla, \D$ are the usual gradient and Laplacian in $\R^n$. Denote $\rho = \frac{1}{2}\ln\frac{1+r}{1-r}$ the hyperbolic distance from the origin to $x\in \H^n$. For hyperbolic radial function, i.e. function $\varphi(\rho)$, there holds
$$\Delta_\H \varphi =\frac{\partial^2 \varphi}{\partial \rho^2} + (n-1)\coth \rho\frac{\partial \varphi}{\partial \rho}.$$

For example, let $f=(\ln\tanh \frac{\rho}{2})^a(\sinh \rho)^b$ with $a, b\in\R$, we have
\begin{align*}
\frac{\Delta_\H f}{f}=b(n-1+b)+\frac{b(n-2+b)}{(\sinh\rho)^2}+\frac{a(a-1)}{(\sinh\rho)^2(\ln\tanh\frac{\rho}{2})^2}+\frac{a(n-2+2b)\coth\rho}{\sinh \rho\ln\tanh\frac{\rho}{2}}.
\end{align*}
Taking $a=\frac{1}{2}$ and $b=\frac{2-n}{2}$, the equality \eqref{E3} means that for any $u \in C_c^1(\H^n \backslash\{0\})$,
\begin{align*}
\int_{\mathbb{H}^n}|\nabla_\H u|^2d v_\H & \geq \frac{n(n-2)}{4}\int_{\mathbb{H}^n} u^2d v_ \H + \frac{(n-2)^2}{4}\int_{\mathbb{H}^n}\frac{ u^2}{(\sinh \rho)^2}d v_\H \\
& \quad + \frac{1}{4}\int_{\mathbb{H}^n}\frac{u^2}{(\sinh\rho)^2(\ln\tanh\frac{\rho}{2})^2}d v_\H.
\end{align*}
This is the optimal Hardy-type inequality showed by \cite[Theorem 2.10]{bggp}.

\medskip
Consider $f=\rho^a (\sinh \rho)^b (\cosh \rho)^\alpha$ with $a, b, \a \in \R$, then
\begin{align}
\label{HyH1}
\begin{split}
\frac{\Delta_\H f}{f} & = \frac{a(a-1)}{\rho^2} + \frac{b(n-2+b)}{\sinh^2\rho} + a(n-1+2b)\frac{\coth\rho}{\rho}\\ & \quad + \alpha (n+2b) + b(n-1+b) + \alpha(\alpha-1)\tanh^2\rho + 2\alpha a \frac{\tanh \rho}{\rho}.
\end{split}
\end{align}

Set first $\alpha=b = 0$ and $a=\frac{2-n}{2}$, applying \eqref{E3}, we have
$$
\int_{\H^n} |\nabla_\H u|^2 dv_\H \geq \frac{(n-2)^2}{4}\int_{\H^n}\frac{u^2}{\rho^2} dv_\H + \frac{(n-1)(n-2)}{2}\int_{\H^n} \frac{\rho\coth\rho - 1}{\rho^2} u^2 dv_\H,
$$
which was given in \cite[Theorem 1.3]{fllm}.

\medskip
Let $a=\frac{1}{2}$ and $b=\frac{1-n}{2}$, \eqref{HyH1} yields
$$
 -\frac{\Delta_\H f}{f} = \frac{1}{4\rho^2}+\frac{(n-1)(n-3)}{4\sinh^2\rho}+\frac{(n-1)^2}{4} + I_\alpha(\rho),
$$
where
$$I_\alpha(\rho):=-\alpha\left[(\alpha-1)\tanh^2\rho+ \frac{\tanh\rho}{\rho}+ 1\right].$$
Taking still $\a = 0$ and using again \eqref{E3}, there holds
\begin{align}\label{2.20}
\begin{split}
\int_{\H^n} |\nabla_\H u|^2 dv_\H & \geq \frac{1}{4}\int_{\H^n}\frac{u^2}{\rho^2}dv_\H + \frac{(n-1)(n-3)}{4}\int_{\H^n}\frac{u^2}{\sinh^2\rho}dv_\H\\
& \quad  + \frac{(n-1)^2}{4}\int_{\H^n} u^2 dv_\H,
\end{split}
\end{align}
which was given by \cite [Theorem 2.1]{bgg} and \cite [Theorem 1.4]{fllm}.

\medskip
Consider now the unit hyperbolic ball ${\rm B}_\H : = \{x \in \B^n,\; \rho < 1\} = \{r < \tanh 1\}$. We want to obtain an improvement of \eqref{2.20} over ${\rm B}_\H$ with a better choice of $\a$.

\medskip
Clearly $I_\alpha(\rho) < 0$ for $\alpha\geq 1$; for $0 < \alpha<1$, we have still $I_\alpha(\rho)<0$ since $\tanh^2\rho < 1$ for any $\rho > 0$. Hence $\a > 0$ are not appropriate. Let $\alpha<0$ and rewrite
$$I_\alpha (\rho)= - \alpha \tanh^2\rho\big[\alpha-1+ J(\rho)\big] \quad \mbox{with } J(\rho):=\frac{\coth\rho }{\rho}+\coth^2\rho.$$
Obviously $J$ is decreasing in $(0,\infty)$, so that for $\rho \in [0, 1]$, $I_\alpha(\rho)\geq G(\a) \tanh^2\rho$ with
$$ G(\alpha) := -\alpha(\alpha - 1) - J(1) \alpha =-\alpha^2 - \big(\coth 1+ {\rm csch}^2 1\big)\alpha.$$
Choosing $\a = \frac{- \coth 1 - {\rm csch}^2 1}{2} < 0$ which reaches the maximum of $G$, we arrive at
\begin{thm}
Let $n \geq 3$. For any $u \in H^1_0({\rm B}_\H)$, there holds,
\begin{align*}
\int_{{\rm B}_\H} |\nabla_\H u|^2 dv_\H & \geq \frac{1}{4}\int_{{\rm B}_\H}\frac{u^2}{\rho^2}dv_\H +  \frac{(n-1)(n-3)}{4}\int_{{\rm B}_\H}\frac{u^2}{\sinh^2\rho}dv_\H \\
& \quad  + \frac{(n-1)^2}{4}\int_{{\rm B}_\H} u^2 dv_\H+\frac{\big(\coth 1 + {\rm csch}^21\big)^2}{4}\int_{{\rm B}_\H}\left(\tanh\rho\right)^2 u^2 dv_\H.
\end{align*}
\end{thm}

\section{Hardy inequality for edge Laplacian}
\reset
Let $X$ be an open set of $\R^n$ and $Y$ be an open set in $\mathbb{R}^q$ containing the origin, $n, q \geq 1$. Consider $\E = (0, 1) \times X \times Y$, which can be regarded as a local model of stretched edge-manifolds
(i.e.~manifolds with edge singularities, see \cite{clw, nsss}) with dimension $N = n+q+1$, and the Riemannian metric
$$t^{-2}dt^2 + dx^2 + t^{-2}dy^2 \quad \mbox{for } (t, x, y) \in \E.$$
Hence the associated gradient and volume form are
$$\nabla_\E = (t\p_t, \p_{x_1}\ldots \p_{x_n}, t\p_{y_1}, \ldots t\p_{y_q}), \quad d\sigma = t^{-1-q}dtdxdy.$$
The corresponding Laplace-Beltrami operator (called edge Laplacian) is the following degenerate elliptic operator:
$$
\Delta_\E = (t\partial_t)^2 -qt\partial_t+\Delta_x +t^2\Delta_y=t^2\partial_{tt} +(1-q) t\partial_t +\Delta_x + t^2\Delta_y.
$$

\begin{thm}
Let $n \geq 2$ and $u\in C_c^1(\E)$, we have
\begin{align}
\label{new-eH}
\int_\E |\nabla_\E u|^2 d\sigma \geq \int_\E \left[\frac{(n-2)^2}{4\psi} + \frac{q^2}{4} + \frac{(n-2)e^2}{8}qW_0\right] u^2d\sigma
\end{align}
where
$$
\psi =e^{-\frac{1}{t^2}}+|x|^2+ |y|^2 \quad \mbox{and} \quad W_0 = \frac{t^{-2}e^{-\frac{1}{t^2}}}{\psi}.
$$
\end{thm}
The above result is motivated by \cite [Proposition 3.5]{clw}, where the following Hardy inequality was used to handle the existence of solution to a Dirichlet problem with edge Laplacian and singular potential.
\begin{align}
\label{eHclw}
\int_\E |\nabla_\E u|^2 d\sigma \geq \frac{n^2}{4}\int_\E W_0 u^2d\sigma, \quad \forall\; u\in C_c^1(\E).
\end{align}
For $n \geq 2$, as $t^{-2}e^{-\frac{1}{t^2}} \leq e^{-1}$ if $t \in (0, 1)$,
\begin{align*}
\frac{W}{W_0} \geq \frac{(n-2)^2e}{4} + \frac{(n-2)e^2}{8}q, \quad \mbox{where  }\; W=\frac{(n-2)^2}{4\psi} + \frac{q^2}{4} + \frac{(n-2)e^2}{8}qW_0.
\end{align*}
Hence \eqref{new-eH} improves greatly \eqref{eHclw}, especially for large $n$ and $q$.

\medskip
\pr Set $f=g(t)\psi^a$ with $a\in \R$, direct computation yields
\begin{align*}
\frac{\Delta_\E f}{f} & = \frac{t^2\partial_{tt}f +(1-q) t\partial_t f +\Delta_x f+ t^2\Delta_y f}{f} \\
& =  t^2 \left[\frac{g''}{g} + 4a\frac{g'}{g\psi t^3}e^{-\frac{1}{t^2}} - \frac{6a}{\psi t^4}e^{-\frac{1}{t^2}}+ \frac{4a}{\psi t^6}e^{-\frac{1}{t^2}} + \frac{4a(a-1)}{\psi^2 t^6}e^{-\frac{2}{t^2}} \right]\\
& \quad + (1-q)t\left[ \frac{g'}{g} + \frac{2a}{\psi t^3}e^{-\frac{1}{t^2}} \right] + \frac{4a(a-1)}{\psi^2}|x|^2 + 2n\frac{a}{\psi}+ t^2\left[ \frac{4a(a-1)}{\psi^2}|y|^2 + \frac{2qa}{\psi}\right]\\
& = \frac{4a(a-1)}{\psi^2}\Big[t^{-4}e^{-\frac{2}{t^2}}+ |x|^2 + t^2|y|^2\Big] + \frac{(2n + 2qt^2)a}{\psi}\\
& \quad + \frac{ae^{-\frac{1}{t^2}}}{t^4\psi}\left[4 + \frac{4{t^3}g'}{g} - 2 (q+2) t^2\right] + \frac{(1 - q)tg'}{g} + t^2\frac{g''}{g}.
\end{align*}
Taking $g(t) = t^b$ with $b\in \R$, there holds
\begin{align*}
\frac{\Delta_\E f}{f} & = \frac{4a(a-1)}{\psi^2}\Big[t^{-4}e^{-\frac{2}{t^2}}+ |x|^2 + t^2|y|^2\Big] + \frac{(2n + 2qt^2)a}{\psi}\\
& \quad  + \frac{ae^{-\frac{1}{t^2}}}{t^2\psi}\left(\frac{4}{t^2}+4b -2q- 4\right) + b(b-q).
\end{align*}
As $t^{-4}e^{-\frac{1}{t^2}} \leq 4e^{-2} < 1$ for $t \in (0, 1)$,
$$
t^{-4}e^{-\frac{2}{t^2}}+ |x|^2 + t^2|y|^2\leq \psi.
$$
Let now $a \leq 0$ and $b = \frac{q}{2}$, we conclude that
\begin{align*}
-\frac{\Delta_\E f}{f} & \geq \frac{4(1-a)a}{\psi} - \frac{(2n + 2qt^2)a}{\psi} - aW_0\left(\frac{4}{t^2} - 4\right)+ \frac{q^2}{4}\\
& \geq  \frac{4(1-a)a - 2na}{\psi} - 2qaW_0t^4e^\frac{1}{t^2} + \frac{q^2}{4}\\
& \geq  \frac{4(1-a)a - 2na}{\psi} -  \frac{e^2aq}{2}W_0 + \frac{q^2}{4}.
\end{align*}
The proof is completed by choosing $a= \frac{2-n}{4}$. \qed

\medskip
Our approach can also improve \eqref{eHclw} in dimension one. Indeed, for $a = 0$, $b = \frac{q}{2}$ and $n = 1$, the above proof still works and implies
$$
\int_\E |\nabla_\E u|^2 d\sigma \geq \frac{q^2}{4}\int_\E u^2d\sigma.
$$

\section{On Heisenberg group}
\label{Heisenberg}
\reset
Let $n \geq 1$ and $\HH_n$ be the standard Heisenberg group, i.e.~$\R^{2n+1}$ endowed with the following group law: For any $\xi = (x, y, t)$, $\tilde{\xi} = (\tilde x, \tilde y, \tilde t) \in \R^n\times \R^n\times \R$,
$$
\xi\circ\tilde{\xi} = \Big({x}+\tilde{x}, y + \tilde y, t+\tilde{t}+2\sum_{i=1}^n(x_i\tilde{y}_i - y_i\tilde{x}_i)\Big).
$$
We equip $\HH_n$ with the norm
$$
\|\xi\|_{\HH_n}:= \rho = \left(r^4 +t^2\right)^\frac{1}{4}, \quad \mbox{where } r^2 = |(x, y)|^2 = |x|^2 + |y|^2.
$$
Define the vector fields
$$
X_i:=\frac{\partial}{\partial x_i}+ 2y_i \frac{\partial}{\partial t},\quad Y_i:=\frac{\partial}{\partial y_i}-2x_i\frac{\partial}{\partial t}, \quad 1\leq i \leq n.
$$
The associated horizontal gradient and Kohn hypoelliptic Laplacian on $\HH_n$ are respectively denoted by $\nabla_{\HH_n}$ and $\Delta_{\HH_n}$, that is
$$\nabla_{\HH_n} = (X_1,\ldots, X_n, Y_1, \ldots, Y_n)$$
and
$$\Delta_{\HH_n} =\sum_{i=1}^{n} X_i^2 + \sum_{i=1}^n Y_i^2 = \Delta_{(x, y)} + 4r^2 \frac{\p^2}{\p t^2} + 4 \sum_{j=1}^n\left(y_j\frac{\p}{\p x_j}-x_j\frac{\p}{\p y_j}\right)\frac{\p}{\p t}.$$
We denote by $Q:=2n+2$ the homogeneous dimension of $\HH_n.$ Here is a series of sharp and new Hardy inequalities.

\begin{thm}
The following estimates hold true.
\begin{equation}\label{HeiH1}
\int_{\HH_n} |\nabla_{\HH_n} u|^2 d\xi \geq \frac{Q^2}{4} \int_{\HH_n} \frac{r^2}{\rho^4} u^2 d\xi + \int_{\HH_n} \frac{r^2}{t^2} u^2 d\xi, \;\; \forall\; u \in C_c^1(\HH_n\backslash\{t = 0\}).
\end{equation}
\begin{equation}\label{HeiH2}
\int_{\HH_n} |\nabla_{\HH_n} u|^2 d\xi \geq \frac{(Q+2)^2}{4} \int_{\HH_n} \frac{r^2}{\rho^4} u^2 d\xi, \;\; \forall\; u \in C_c^1(\HH_n\backslash\{0\});
\end{equation}
\begin{equation}\label{HeiH4}
\int_{\HH_n} |\nabla_{\HH_n} u|^2 d\xi \geq \frac{(Q-4)^2}{4} \int_{\HH_n} \frac{u^2}{r^2} d\xi + 9 \int_{\HH_n} \frac{r^2}{\rho^4} u^2 d\xi, \;\; \forall\; u \in C_c^1(\HH_n\backslash\{r = 0\}).
\end{equation}
\end{thm}

\medskip
Here $d\xi$ denotes the Lebesgue measure in $\R^{2n+1}$. The inequality \eqref{HeiH1} improves \cite[Theorem 1.1]{ly} where Luan-Yang proved similar estimate with $\frac{Q^2 - 4}{4}$ instead of the sharp coefficient $\frac{Q^2}{4}$.

\medskip
\pr To prove \eqref{HeiH1}, clearly we need only to consider $u$ supported in $\HH_n^+$, the half space where $t > 0$. With integration by parts, we can check easily that \eqref{E3} becomes
$$\int_{\HH_n^+} |\nabla_{\HH_n} u|^2 d\xi  =  \int_{\HH_n^+} \frac{-\Delta_{\HH_n} f}{f}u^2 d\xi + \int_{\HH_n^+} \left|\nabla_{\HH_n} u - \frac{u}{f}\nabla_{\HH_n} f \right|^2 d\xi.$$
Let $f=\rho^a r^b t^c$ with $a$, $b$, $c\in \mathbb{R}$, there holds
$$
-\frac{\Delta_{\HH_n} f}{f} = -\frac{b(Q-4+b)}{r^2} -a(Q-2+2b+4c+a)\frac{r^2}{\rho^4}-4c(c-1)\frac{r^2}{t^2}.
$$
Set $a=-\frac{Q}{2}$, $b=0$ and $c=\frac{1}{2}$, we get \eqref{HeiH1}.

\medskip
The estimates \eqref{HeiH2} and \eqref{HeiH4} are obtained similarly, with different choice of parameters in $f$, working on the corresponding domain.
\begin{itemize}
\item \eqref{HeiH2} is obtained with $a=-\frac{Q+2}{2}$, $b=0$ and $c=1$.
\item Taking $a=-3$, $b=-\frac{Q-4}{2}$ and $c=1$, \eqref{HeiH4} is proved. \qed
\end{itemize}

\section{Gurshin type operators}
\reset
Let $\xi =(x, y)\in \mathbb{R}^{d+k}$ with $d, k\geq 1$. Let $\gamma > 0$, consider
$$
X_i:=\frac{\partial}{\partial x_i},\; Y_j:=|x|^\gamma \frac{\partial}{\partial y_j}, \quad 1 \leq i \leq d, \; 1 \leq j \leq k.
$$
The associated Grushin gradient and Grushin Laplacian are defined by
$$
\nabla_\gamma =(X_1,...,X_d,Y_1,...,Y_k) = (\nabla_x, |x|^\gamma\nabla_y), \quad \Delta_\gamma =\sum_{i=1}^d X_i^2 +\sum_{j=1}^k Y_j^2 = \Delta_x+|x|^{2\gamma}\Delta_y.
$$
Here we denote by $\rho$ the Grushin distance from the origin to $\xi = (x, y)$:
$$\rho = \|\xi\|_G = \Big[|x|^{2+2\gamma} + (1+\gamma)^2|y|^2\Big]^{\frac{1}{2+2\gamma}}.$$
D'Ambrosio in \cite [Theorem 3.3]{am2} proved a well-known Hardy inequality (among many others) for Grushin operators: For any $u \in C_c^1((\R^d\backslash\{0\})\times \R^k)$,
\begin{align}
\label{GHA}
\int_{\R^{d+k}}|\nabla_\gamma u|^2 d\xi \geq \frac{(d-2)^2}{4}\int_{\mathbb{R}^{d+k}}\frac{u^2}{|x|^2}d\xi \geq \frac{(d-2)^2}{4}\int_{\mathbb{R}^{d+k}}\frac{u^2}{\rho^2}d\xi.
\end{align}
When $d\geq 3$, applying density argument, the above estimate holds true for $u \in D^{1,2}_\gamma(\mathbb{R}^{d+k})$, the closure of $C_c^1(\mathbb{R}^{d+k})$ endowed with the norm $\|\nabla_\gamma u\|_{L^2}$. In the following, we show some examples of improved Hardy inequalities.

\begin{thm}
For any $u \in C_c^1\left((\R^d\backslash\{0\})\times (\R^k\backslash\{0\})\right)$, we have
\begin{align}\label{GH1}
\begin{split}
\int_{\R^{d+k}}|\nabla_\gamma u|^2 d\xi & \geq \frac{(d-2)^2}{4}\int_{\R^{d+k}}\frac{u^2}{|x|^2}d\xi + \frac{(k-2)^2}{4}\int_{\R^{d+k}}\frac{|x|^{2\gamma}}{|y|^2}u^2 d\xi\\
& \quad + (1+\gamma)^2 \int_{\R^{d+k}}\frac{|x|^{2\gamma}}{\rho^{2+2\gamma}} u^2 d\xi;
\end{split}
\end{align}
\begin{equation}\label{GH2}
\int_{\R^{d+k}}|\nabla_\gamma u|^2 d\xi\geq \frac{(d-2)^2}{4}\int_{\R^{d+k}}\frac{u^2}{|x|^2}d\xi + \frac{k^2(1+\gamma)^2}{4} \int_{\R^{d+k}}\frac{|x|^{2\gamma}}{\rho^{2+2\gamma}} u^2 d\xi;
\end{equation}
\begin{equation}\label{GH3}
\int_{\R^{d+k}}|\nabla_\gamma u|^2 d\xi\geq \frac{(k-2)^2}{4}\int_{\R^{d+k}}\frac{|x|^{2\gamma}}{|y|^2}u^2 d\xi + \frac{(d+2\gamma)^2}{4} \int_{\R^{d+k}}\frac{|x|^{2\gamma}}{\rho^{2+2\gamma}} u^2 d\xi.
\end{equation}
\end{thm}

Clearly, \eqref{GH2} is valid for $u \in C_c^1((\R^d\backslash\{0\})\times \R^k)$, hence improves \eqref{GHA}. When $d, k \geq 3$, \eqref{GH1}--\eqref{GH3} hold true for $u \in D^{1,2}_\gamma(\mathbb{R}^{d+k})$, in particular \eqref{GH1} improves also \eqref{GHA} if $k \geq 3$.

\medskip
\pr Let $f= |x|^a |y|^b \rho^c$ with $a, b, c \in \R$, there holds
\begin{align*}
\frac{\Delta_\gamma f}{f} & = \frac{\Delta_x f+ |x|^{2\gamma}\Delta_y f}{f}\\& = \frac{a(d-2+a)}{|x|^2} + \frac{b(k-2+b)|x|^{2\gamma}}{|y|^2} + c\big[c-2+d+2a+(1+\gamma)(k+2b)\big]\frac{|x|^{2\gamma}}{\rho^{2+2\gamma}}.
\end{align*}
Using integration by parts, we can apply \eqref{E3} with $\nabla_\gamma$, $\Delta_\gamma$ and $d\xi$.
\begin{itemize}
\item Choosing $a=\frac{2-d}{2}$, $b=\frac{2-k}{2}$ and $c=-(1+\gamma)$, we get \eqref{GH1}.
\item \eqref{GH2} is given by $a =\frac{2-d}{2}$, $b=0$ and $c=-\frac{k(1+\gamma)}{2}$.
\item We obtain \eqref{GH3} with $a=0,~b=\frac{2-k}{2}$ and $c=-\frac{d+2\gamma}{2}$. \qed
\end{itemize}

\section{Further remarks}
\reset

\subsection{Bessel pair and capacity}
By \eqref{E2}, we can generalize Ghoussoub-Moradifam's definition of Bessel pair: $(V,W)$ is called a Bessel pair over an open set $\Omega \subset ({\mathcal M}, g)$ if the following holds true:
$$\mbox{There exists a positive function $f\in C^1(\Omega)$ such that $-{\rm div}_g(V\nabla_g f)\geq W f$ in $\Omega$.}$$
Applying the equality \eqref{E2}, we see that
$$
\int_{\Omega} V\|\nabla_g u\|^2 dg \geq \int_{\Omega} W u^2 dg,\quad \quad \forall~u\in C_c^1(\Omega).
$$

When $V, W\in C^1(\Omega\backslash\Sigma)\cap L_{loc}^1(\Omega)$ forms a Bessel pair over $\Omega\backslash \Sigma$ with a negligible set $\Sigma \subset\subset \Omega$, a natural question is to know when we can extend the couple $(V, W)$ as a Bessel pair over $\Omega$. Here is a sufficient condition in the spirit of \cite{ma1}.
\begin{defn}
Let $\Omega \subset ({\mathcal M}, g)$ be open and $\Sigma \subset\subset \O$. We say that ${\rm cap}_V(\Sigma, \O)=0$ if ${\rm vol}_g(\Sigma) = 0$, and for any open set $\Omega'$ satisfying $\Sigma\subset\subset\Omega'\subset\Omega$, there holds ${\rm cap}_{V, \Omega'}(\Sigma)=0$ where
$$
{\rm cap}_{V, \Omega'}(\Sigma)=\inf\left\{\int_{\Omega'} V\|\nabla_g \eta\|^2 dg, \; \eta \in C_c^1(\Omega'), \; \Sigma \subset\subset\{\eta = 1\} \right\}.
$$
\end{defn}

The following result could be well-known, we show a proof here for the convenience of readers.
\begin{prop}
Let $V, W \in C^1(\Omega\backslash \Sigma)\cap L_{loc}^1(\Omega)$ be a Bessel pair over $\Omega\backslash\Sigma$ where $\Sigma \subset\subset \Omega$. If ${\rm cap}_V(\Sigma, \O) =0$, then $(V, W)$ is a Bessel pair in $\O$.
\end{prop}

\pr Given $\delta>0$, let $\Omega_{\delta}={\{x\in \Omega,~ d(x,\Sigma)<\delta}\}$, then $\Omega_\delta\subset \Omega$. As ${\rm cap}_{V, \Omega_\delta}(\Sigma)=0$, there exists a sequence $(\eta_k) \in C_c^1(\Omega_\delta)$ such that
$$\Sigma \subset\subset\{\eta_k = 1\}, \quad \lim_{k\to\infty}\int_{\Omega_\delta} V\|\nabla_g \eta_k\|^2 dg = 0.$$
For $u\in C_c^1(\Omega)$, $u\zeta_k\in C_c^1(\Omega\backslash\Sigma)$ with $\zeta_k = 1 - \eta_k$. Therefore
\begin{equation}
\int_{\Omega\backslash\Sigma} V \|\nabla_g(u\zeta_k)\|^2 dg \geq \int_{\Omega\backslash\Sigma} W(u\zeta_k)^2 dg.
\end{equation}
Moreover, for any $\e > 0$, there exists $C_\e > 0$ such that
\begin{align*}
\int_{\Omega\backslash\Sigma} V \|\nabla_g(u\zeta_k)\|^2 dg &\leq (1+\epsilon)\int_{\Omega} V \|\nabla_g u\|^2 dg + C_\epsilon \int_{\Omega} V u^2\|\nabla_g \zeta_k\|^2 dg \\
& \leq (1+\epsilon)\int_{\Omega} V \|\nabla_g u\|^2 dg + C_\epsilon \|u\|^2_\infty\int_{\Omega} V \|\nabla_g \eta_k\|^2 dg.
\end{align*}
Taking $k\to \infty$, there holds
$$
(1+\epsilon)\int_{\Omega} V \|\nabla_g u\|^2 dg \geq \int_{\Omega\backslash\Sigma} W(u\zeta_k)^2 dg \geq \int_{\Omega\backslash\Omega_\delta} Wu^2 dg.
$$
On the other hand,
$$
\lim_{\delta\to 0^+}\int_{\Omega\backslash\Omega_\delta} Wu^2 dg = \int_{\Omega\backslash\Sigma} W u^2 dg = \int_{\Omega} W u^2 dg.
$$
Combining the above estimates and let $\delta \to 0^+$,
$$
(1+\epsilon)\int_{\Omega} V\|\nabla_g u\|^2 dg \geq \int_{\Omega} W u^2 dg.
$$
Taking $\epsilon \to 0^+$, we are done. \qed

\subsection{$L^p$ Hardy inequalities}
Until now, we concentrate our discussion for Hardy inequalities in $L^2$ setting, here we point out that the idea to use equalities works also for general $L^p$ setting with $p > 1$.

\medskip
Let $u \in C_c^1(\O)$, $p > 1$, $V\in C^1(\Omega)$ be a nonnegative weight. Consider $\vec{F}\in C^1(\Omega, T_g{\mathcal M})$ and a family of inner product $\langle \cdot, \cdot\rangle \in C^1(\O, \Lambda^2T_g{\mathcal M})$.
\begin{align}
\label{Ep}
\begin{split}
& \quad \int_{\Omega} V\|\nabla_g u\|^p dg\\
& = - \int_{\Omega} V\left[(p-1)\|\vec F\|^p |u|^p + \|\vec F\|^{p-2} \langle\vec F, \nabla_g (|u|^p)\rangle\right] dg + \int_{\Omega} V{\mathcal R}(u, \vec F)dg \\
 & = \int_{\Omega} \left[{\rm div}_g\left(V \|\vec F\|^{p-2}\vec F\right)-(p-1)V\|\vec F\|^p\right] |u|^p dg + \int_{\Omega} V{\mathcal R}(u, \vec F)dg,
\end{split}
\end{align}
where
\begin{align*}
{\mathcal R}(u, \vec F) & = (p-1)\|\vec F\|^p |u|^p + \|\vec F\|^{p-2} \langle\vec F, \nabla_g (|u|^p)\rangle + \|\nabla_g u\|^p\\
& = (p-1)\|\vec F\|^p |u|^p + \|\nabla_g u\|^p + p\|\vec F\|^{p-2} |u|^{p-2}u \langle\vec F, \nabla_g u \rangle \geq 0.
\end{align*}
The last estimate is given by the Cauchy-Schwarz inequality for inner product and Young's inequality. Therefore, we obtain a $L^p$-Hardy inequality
\begin{align}
\label{HLp}
\int_{\Omega} V\|\nabla_g u\|^p dg \geq \int_{\Omega} W|u|^p dg,
\end{align}
with
\begin{align*}
W = {\rm div}_g\left(V \|\vec F\|^{p-2}\vec F\right)-(p-1)V\|\vec F\|^p.
\end{align*}
In particular, if $\vec F = -\frac{\nabla_g f}{f}$ with $f>0$, there holds
\begin{align}
\label{WLp}
W = -\frac{{\rm div}_g\left(V\|\nabla_g f\|^{p-2}\nabla_g f\right)}{f^{p-1}} =:-\frac{L_{V,p}f}{f^{p-1}}.
\end{align}

The following examples are direct consequence of the estimate \eqref{HLp}, hence somehow of the equality \eqref{Ep}.
\begin{itemize}
\item Let $V = 1$, $p > 1$ and $f=|x|^{1 - \frac{n}{p}}$ in $\R^n$, we have
$$-\frac{{\rm div}(|\nabla f|^{p-2} \nabla f)}{f^{p-1}} = \left|\frac{n-p}{p}\right|^p |x|^{-p} \quad \mbox{in } \R^n\backslash\{0\},$$
which means that
$$
\int_{\R^n} |\nabla u|^p dx \geq \left|\frac{n-p}{p}\right|^p \int_{\R^n} \frac{u^p}{|x|^p} dx,\quad \forall~u\in C_c^1(\R^n\backslash\{0\}).
$$
This is the well-known $L^p$-Hardy inequality, which generalizes \eqref{Hardy0} to any dimension.
\item Consider $\langle \vec v, \vec w \rangle_A = \langle A(x)\vec v, \vec w\rangle$ in $\R^n$ with a smooth family of positive definite symmetric matrix $A(x)$ and $V = 1$. Now $$L_{V,p}f = {\rm div}(A(x)|\nabla f|^{p-2}\nabla f) =: \L_{A,p}f.$$
Take $f = E^{1-\frac{1}{p}}$ with a positive function $E$ and $-\L_{A,p}E = \mu$, we derive from \eqref{HLp}--\eqref{WLp} that for any $u \in W^{1, p}_0(\O)$, $p > 1$,
\begin{align*}
\int_\O |\nabla u|_A^p dx \geq \left(\frac{p-1}{p}\right)^p\int_\O \frac{|\nabla E|_A^p}{E^p}|u|^p dx + \left(\frac{p-1}{p}\right)^{p-1}\int_\O \frac{|u|^p}{E^{p-1}} d\mu
\end{align*}
which is the inequality (9.34) in \cite{gm2}. If moreover $\mu$ is a nonnegative finite measure on $\O$, we get
$$
\int_\O |\nabla u|_A^p dx \geq \left(\frac{p-1}{p}\right)^p\int_\O \frac{|\nabla E|_A^p}{E^p}|u|^p dx,  \quad \forall\; u \in C_c^1(\O).
$$
which is given in \cite{as}.
\item Let $\HH_n$ be the Heisenberg group. Here all the notations are that in section \ref{Heisenberg}. Let $V = r^{\beta-p}\rho^{2p-\alpha}$ with $\alpha, \beta\in \mathbb{R}$, $p > 1$, $n\geq 1$ and $f=\rho^b$. As $|\nabla_{\HH_n}\rho| = \frac{r}{\rho}$ for $\rho > 0$, there holds, in $\HH_n\backslash\{0\}$,
\begin{align*}
-\frac{{\rm div_{\HH_n}}(V|\nabla_{\HH_n} f|^{p-2} \nabla_{\HH_n} f)}{f^{p-1}} = -|b|^{p-2}b\Big[2n+2+b(p-1)+\beta-\alpha\Big] \frac{r^\beta}{\rho^\alpha}.\end{align*}
Choosing $b=-\frac{2n+2+\beta-\alpha}{p}$, we arrive at
$$
\int_{\HH_n}r^{\beta-p}\rho^{2p-\alpha} |\nabla_{\HH_n} u|^p dx \geq \left(\frac{2n+2+\beta-\alpha}{p}\right)^p \int_{\HH_n} \frac{r^\beta}{\rho^\alpha}u^p dx
$$
for any $u\in C_c^1(\HH_n\backslash\{0\})$. This enables us Theorem 3.2 in \cite{am1}.
\end{itemize}

\bigskip
\noindent{\bf Acknowledgements.} When the second author was an undergraduate student, Professor Chaojiang Xu taught him the theory of distribution and partial differential equations, he would like to express his gratitude to Professor Xu for the guidance and constant support. Both authors are partially supported by Science and Technology Commission of Shanghai Municipality (STCSM), grant No. 18dz2271000. X. H. is partially supported by NSFC (No. 11971169). The authors thank also the anonymous referees for their useful comments.

\end{document}